\documentclass[11pt]{article} 
\usepackage{graphicx,color}
\usepackage{amsthm,amssymb,amsfonts, amsmath}
\usepackage{varioref}
\usepackage[bookmarksnumbered,plainpages]{hyperref} 
\usepackage{epsfig, url}

\newtheorem{theorem}{Theorem}[section]

\newtheorem{lemma}{Lemma}[section]

\newtheorem{definition}{Definition}[section]

\newtheorem{rem}{Remark}[section]

\DeclareMathOperator{\tg}{tan}

\begin{document}
\hoffset=-1.2cm
\voffset=-2.5cm

\def\Om{\Omega}
\def\om{\omega}
\def\e{\varepsilon}
\def\all{\alpha^{\e}}
\def\al#1{\alpha_{#1}^{\e}}
\def\bt#1{\beta_{#1}^{\e}}
\def\lm{\lambda}
\def\lme{\lambda_\e}
\def\me{\mu_\e}
\def\intl{\int\limits}
\def\btt{\beta^{\e}}
\def\f{\varphi}
\def\Oe{\Om^{\e}}
\renewcommand\le{\leqslant}
\renewcommand\leq{\leqslant}
\renewcommand\geq{\geqslant}
\renewcommand\ge{\geqslant}
\renewcommand\kappa{\varkappa}

\title
{Complete WKB asymptotics of high frequency vibrations\\ %
in a stiff problem\footnote{ %
Published in \emph{Mat. Stud. 14, no.1 (2001): 59--72} %
}
\footnote{ %
English is improved in 2008
}}

\author{Natalia Babych and Yuri Golovaty}

\date{}

\maketitle

\begin{abstract}
Asymptotic behaviour of eigenvalues and eigenfunctions
of a stiff problem is described 
in the case of the fourth-order ordinary differential operator. 
Considering the stiffness coefficient that depends on a small
parameter $\e$ and vanishes as $\e\to 0$ on a subinterval,
we prove the existence of low and high frequency resonance vibrations.
The low frequency vibrations admit the power series expansions on $\e$
but this method is not applicable to the description of high frequency vibrations.
However, the non-classical asymptotics on $\e$ of the high frequency vibrations
were constructed using the WKB method.
\end{abstract}

\textbf{MSC}: Primary 34E20; Secondary 74K10

\section*{Introduction and main results.}
Stiff vibrating systems belong to a class of systems with singularly perturbed
potential energy. Stiff problems are known in particular as
boundary value problems for
 differential equations with very contrasting
values of coefficients in different sub-domains.
They relate to modelling vibrations of  elastic
systems consisting of two (or more) materials with one of them being very stiff with respect to the other. 
For the first time the stiff problems were 
investigated by J.~L.~Lions~[1].

However, it is also of interest to describe the asymptotic behaviour
of spectral properties for the stiff problems. 
These system has two types of
eigenvibrations, namely low frequency vibrations and
high frequency ones. From a physical viewpoint we can postulate that two
kinds of eigenvibrations can appear: one for the stiffer structure
and the other for the softer structure.
Different aspects of the spectral stiff problems are considered in [2]--[10]
with the best general reference being [8].
The asymptotic behaviour of the low frequency vibrations has been widely
studied with different techniques [2], [3], [7] and [9].

In this paper, following [10] we consider the phenomenon of high frequency
vibrations. The leading terms of high frequency vibrations
for different problems were constructed in [7-9]. Information on the behaviour
of high frequency vibrations was also provided in [10].
Studying the stiff problem for the forth-order differential operator,
we construct the complete asymptotic expansions of
high frequency vibrations using the WKB technique.
 
\vspace{5mm}
 
Why do high frequency vibrations appear? On the one hand,
the spectrum $\{\lm_i^{\e}\}_{i=1}^\infty$ of the stiff problem
is asymptotically dense in $[0,\infty)$ as $\e\to 0$ (see [10]).
On the other hand, for fixed $i$ the asymptotic expansions of an eigenvalue
$\lm_i^{\e}$  and the corresponding eigenfunction $u_i^{\e}$ are nonuniform
with respect to $i$. 
Hence the sequences $u^\e_{i(\e)}$ with
$i(\e)\to \infty$ can support stable forms of vibration as $\e\to 0$.
In addition, the corresponding sequencies of eigenvalues $\lm_i^{\e}$ converge
to some positive limit points (see Fig. 1, 2). It is also shown that
the approximation of the limit form of vibrations
by sequence $\{u_i^{\e}\}$ has the discrete character.
Therefore we construct and
justify the asymptotics for a family of 
discrete sets of a small parameter.

\vspace{5mm}

\hspace{5mm}
\includegraphics[width=.9\textwidth]{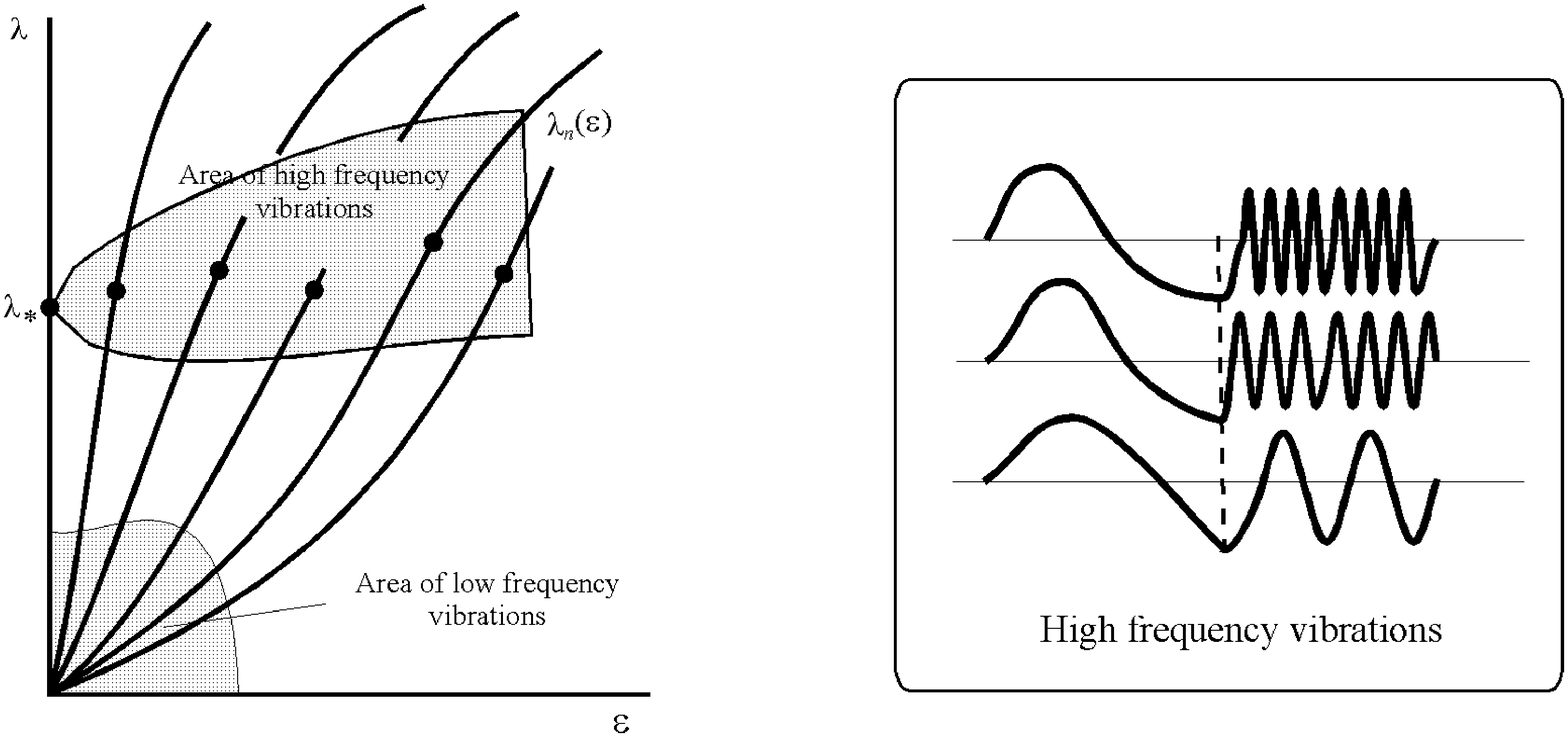}

\hspace{4cm}
Fig. 1 \hspace{4cm} Fig. 2

\vspace{5mm}

Moreover, such approximation is  ambiguously determined.
Then we introduce a deformation parameter $\delta$ in the asymptotics.
Consequently, for each  $\delta$  we  construct
expansions of $u^\e$, though they are asymptotically equivalent when $\e\to 0$.

\section{Problem statement}
Let an interval $(a,b) \subset \mathbb R$ contains the origin.
Let us consider the eigenvalue problem
\begin{eqnarray}
 &
 \displaystyle
 \frac{d^2}{dx^2}\left( k_\e(x)\frac{d^2u_\e}{dx^2}\right)
 - \lme u_\e(x) = 0,
      \quad x \in (a,b),
 &
 \label{tag1}\\
 &
 u_\e(a) = u'_\e(a) = 0,
 \quad 
 u_\e(b) = u'_\e(b)=0,
 &
 \label{tag2}\\[2mm]
 &
 u_\e(-0) = u_\e(+0),
 \quad u'_\e(-0) = u'_\e(+0),
 &
 \label{tag3}\\
 &
 (k_0(x) u''_\e)(-0) = \e^4 (k_1(x) u''_\e)(+0),\quad
 (k_0(x) u''_\e)'(-0) = \e^4 (k_1(x) u''_\e)'(+0). 
 &
 \label{tag4}
\end{eqnarray}
Here $\e$ is a small positive parameter and a function
$$
  k_\e(x)=
  \left\{ 
  \begin{array}{rl}
  k_0(x),& a < x <0 
  \\[1mm]
  \e^4 k_1(x),& 0 < x < b
  \end{array}
  \right.
$$
is  smooth and strictly positive in $(a,0)\cup (0,b)$. We study  the asymptotic behaviour of
the eigenvalues $\lme$ and the eigenfunctions $u_\e$ 
of (1)--(4) as $\e \to 0$.

\vspace{5mm}

The problem models eigenvibrations of a non-homogeneous rod.
The rod consists of two components with the same density of mass and
strongly  different elastic properties. The fourth power of $\e$ in
the definition of
$k_\e$ is suitable for the next consideration (see Section 4).

\vspace{5mm}

Let us introduce the Sobolev space $H^2_0(a,b)$ as the closure
of set $C_0^{\infty}(a,b)$ with respect to the norm
$$
\|u\|=\left(\int_a^b |u''|^2\,dx\right)^{1/2},
$$
and bilinear forms in $H^2_0(a,b)$
$$
  a_0(\f,\psi)=\int_a^0 k_0 {\f}''{\psi}''\,dx,\quad
  a_1(\f,\psi)=\int_0^b k_1 {\f}''{\psi}''\,dx,\qquad
  a_\e=a_0+\e^4a_1.
$$
For each $\e>0$ we note by  $\|\cdot\|_\e$ the norm
in $H^2_0(a,b)$ associated with the form $a_\e(\cdot,\cdot)$.
It is known as an energetic norm associated with the elastic-like energy.

Let us consider a variational formulation of (1)--(4) that is
to find 
$\lm_\e$ and $u_\e\in H^2_0(a,b)$ such that
\begin{equation}
  a_\e(u_\e,\f)-\lm_\e(u_\e,\f)_{L_2(a,b)}=0,
  \qquad
  \f\in H^2_0(a,b).
  \label{tag5}
\end{equation}
Problem (5) is a standard eigenvalue problem with a real discrete spectrum.
For each fixed $\e$ let us consider an eigenvalue sequence
$$
 0 < \lm^\e_1 < \lm^\e_2 < \dots < \lm^\e_i < \dots,
 \quad\text{where}\quad
 \lm^\e_i\to\infty
 \quad\text{as}\quad 
 i\to\infty.
$$
Note that each eigenvalue is simple.
Let the corresponding eigenfunctions
$\{u^\e_i\}_{i=1}^{\infty}$ form 
an orthonormal basis in $L_2(a,b)$.

\section{Asymptotic expansions of low frequency vibrations}
We investigate the asymptotic behaviour of  eigenvalues $\lm^\e_i$
and eigenfunctions $u^\e_i$ for a fixed number $i$.

\begin{lemma}
Each eigenvalue $\lm^\e_i$ is a continuous
function with respect to the parameter $\e$, $\e\in(0,1)$.
Moreover,
$$
\lm^\e_i\leq C_i\e^4,
$$
where constant $C_i$ is independent of $\e$.
\end{lemma}

\proof
The continuity of eigenvalues follows from the variational principle
\begin{equation}
 \lm^\e_i=\sup_P\inf_{f\in P^{\perp}\setminus\{0\}}
 \frac{a_0(f,f)+\e^4 a_1(f,f)}{\|f\|^2_{L_2(a,b)}},
	\label{tag6}
\end{equation}
where $P$ is a $(i-1)$-dimensional subspace of $H^2_0(a,b)$ and $P^{\perp}$
is the orthonormal complement of P.

Suppose that the supremum in (6) 
is achieved on a subspace $P_i$.
Since $P_i$ is finite dimensional, we choose a function
$f_i\in P^{\perp}_i\setminus\{0\}$ that vanishes in $(a,0)$.
Then
\begin{align*}
 \lm^\e_i & = \inf_{f\in P^{\perp}_i\setminus\{0\}}
 \frac{a_0(f,f)+\e^4 a_1(f,f)}{\|f\|^2_{L_2(a,b)}}
 \\
 & \leq
 \frac{a_0(f_i,f_i)+\e^4 a_1(f_i,f_i)}{\|f_i\|^2_{L_2(a,b)}}=
 \e^4\frac{a_1(f_i,f_i)}{\|f_i\|^2_{L_2(0,b)}}=C_i\e^4,
\end{align*}
since $a_1(f_i,f_i)\not =0$.
\hfill
$\Box$

\vspace{5mm}

We postulate the expansions of the eigenvalue $\lme=\lambda_i^\e$ and the
eigenfunction $u_{\e}=u^{\e}_i$ for given $i \in \mathbb N$:
\begin{equation}
\begin{array}{r l l}
 \lme &\sim& \e^4(\lm_0+\e^4\lm_1+\dots),
 \\[2mm]
 u_{\e}(x)&\sim&\left\{
 \begin{array}{ll}
              \e^4(u_0(x)+\e^4 u_1(x)+\dots),&\quad x\in(a,0),\\
              v_0(x)+\e^4 v_1(x)+\dots,&\quad x\in(0,b).
 \end{array}
 \right.
 \end{array}
 \label{tag7}
\end{equation} 
Substituting (7) into (1)--(4) we deduce that $\lm_0$ is an eigenvalue and $v_0$
is an eigenfunction of the problem
\begin{equation}
\begin{array}{c}
\displaystyle
 (k_1 v''_0)''(x)-\lm_0 v_0(x)=0,\quad x\in(0,b),
 \\[2mm]
\displaystyle
 v_0(0)=v'_0(0)=0,\quad v_0(b)=v'_0(b)=0.
\end{array}
\label{tag8}
\end{equation}
Note that all eigenvalues of the problem are simple.

Then the function $u_0$ is a solution of the boundary value problem
$$
\begin{array}{c} 
\displaystyle
(k_0 u''_0)''(x) = 0,\quad x \in (a,0),
\\[2mm]
\displaystyle
 u_0(a) = u'_0(a) = 0,\quad k_0(0) u''_0(0) = k_1(0) v''_0(0),
        \quad (k_0 u''_0)'(0)=(k_1v''_0)'(0).
\end{array}
$$
Let us define next terms of expansions (7).
The function $v_1$ satisfies
$$
\begin{array}{c}
\displaystyle
 (k_1 v''_1)''(x) - \lm_0 v_1(x) = \lm_1 v_0(x),
 \quad x\in(0,b),
 \\[2mm]
\displaystyle
 v_1(0)=u_0(0),
 \quad v'_1(0) = u_0(0),
 \quad v_1(b) = v'_1(b)=0.
\end{array}
$$
Since $\lm_0$ is an eigenvalue of (8), such problem has
a solution under the condition
$$
 \lm_1
 =
 \| v_0 \|_{L_2(0,b)}^{-1} ((k_1 v_0'')' u_0 - k_1 v_0'' u_0') \Big|_{x=0}.
$$
Further, the function $u_1$ is a solution of
$$
\begin{array}{c}
\displaystyle
(k_0 u''_1)''(x) = \lm_0 u_0(x),
\quad x\in(a,0),
\\[2mm]
\displaystyle
 u_1(a)=u'_1(a)=0,
 \quad k_0(0) u''_1(0) = k_1(0) v''_1(0),
        \quad (k_0 u''_1)'(0) = (k_1 v''_1)'(0).
\end{array}
$$

The general terms of (7) can be found in the same way. The asymptotic
expansions are justified by classical methods [14, 2].

\vspace{5mm}

The vibrations described above are corresponding
to the low levels of potential energy, since $a_\e(u_\e,u_\e)=o(\e^2)$ as $\e\to 0$.
Naturally we have vibrations of the soft part of a system only. The soft part
is clamped at $x=0$ by the immovable stiffer structure. Note that
the leading terms of
expansions (7) are determined by the soft part.
Since every eigenfunction $u^\e_i$  converges to an eigenfunction $v_0$ of (8)
extended by zero to $(a,b)$, 
the set $\{u^\e_i\}_{i=1}^{\infty}$ is not a basis in
$L_2(a,b)$ for $\e=0$.  In other words, the ``low frequency" region does not
provide a good insight on the vibration problem over all $(a,b)$. Therefore
we consider here another kind of vibrations, namely, vibrations with 
"finite non-vanishing energy".

\section{High frequency vibrations}
On the one hand the asymptotics of $\lm^\e_i$ and $u^\e_i$ are nonuniform with
respect to $i$, on the other hand the spectrum of (1)--(4) is asymptotically
dense in the positive spectral semi-axis (see Lemma 2).
On account of the above remarks, we can find stabile vibrations
$u_{i(\e)}^{\e}$  associated with certain sequences of eigenvalues $\lm_{i(\e)}^{\e}$ with
$i(\e)\to\infty$.

Let us denote by $E$ the set of all pairs $(\e,\lm(\e))$, where $\lm(\e)$
is an eigenvalue of (1)--(4) for some $\e\in (0,1)$.

\begin{lemma} 
The closure of the set $E$ includes semi-axis
$\{(\e,\lm): \e=0, \lm\ge0\}$. That is to say each positive  $\lm$
can be approximated by a sequence of eigenvalues  $\lm(\e)$.
\end{lemma}

\proof
Suppose that for a certain  $\lm>0$
there exists a neighborhood $U$ of the point $(0,\lm)$ 
such that $U\cap E=\varnothing$. 
Let us choose  $\e^*>0$ small enough for $(\e^*,\lm)\in U$.
Recall that  $\lm_i(\e)$ is a continuous function with respect to $\e$ and
$\lm_i(\e)\to 0$ as $\e\to 0$ (see Fig. 1). 
Since the neighborhood $U$ by definition does not contain points
$(\e^*,\lm_i(\e^*))$
then $\lm_i(\e^*)\leq\lm$ for all $i\in \mathbb N$
that is impossible.
\hfill $\Box$

\vspace{5mm}

Let us consider a sequence of pairs $(\e_i,\lm_i)\in E$ 
that converges in
$\mathbb R^2$ to $(0,\om^4)$ with $\om>0$ as  $i\to\infty$. 
Let $u_i$ be an eigenfunction of (1)--(4), being
normalized in $L_2(a,b)$ and associated with $\lm_i$.

\begin{definition}
We say that the eigenfunction sequence $\{u_i\}_{i=1}^{\infty}$ supports
high frequency vibrations as $i\to\infty$, if the sequence
of restrictions $u_i|_{(a,0)}$ has a nonzero limit $v$ in $H^2(a,0)$.
Then $\om$ 
is referred to as the limit frequency and 
$v$ as the limit form of these vibrations.
\end{definition}

\begin{rem} 
It will be shown that the functions $u_i$ have
strongly oscillatory character in $(0,b)$. Moreover, the vibrations $u_i$ have the ``energy" close to $\om^4$ as $i\to\infty$.
\end{rem}

\begin{theorem}
If a sequence $\{u_i\}_{i=1}^{\infty}$ supports 
high frequency vibrations with a limit frequency $\om$ and
a limit form $v$, then $\om^4$ is an eigenvalue and
$v$ is an eigenfunction of the problem
\begin{equation}
 \gathered
 \left( k_0(x)v''\right)''-\om^4 v=0,\quad
   x\in (a,0),\\
 v(a)=v'(a)=0,\quad v''(0)=v'''(0)=0.
 \endgathered
 \label{tag9}
\end{equation}
and restrictions $u_i|_{(0,b)}$ converge to zero  in the weak
topology of $L_2(0,b)$.
\end{theorem}

\proof
From (5) we have $\|u_i\|_\e=\lm_i\|u_i\|_{L_2(a,b)}=\lm_i$. Then
\begin{equation}
\e_i^2a_1(u_i,\f)\leq\e_i^2a_1^{1/2}(u_i,u_i)a_1^{1/2}(\f,\f)
\leq \|u_i\|_{\e_i} a_1^{1/2}(\f,\f)\le\lm_i a_1^{1/2}(\f,\f)
\label{tag10}
\end{equation}
for any $\f\in H^2_0(a,b)$.
Let us consider identity (5) for test functions $\f\in C_0^{\infty}(0,b)$
$$
 (u_i,\f)_{L_2(0,b)}=\e^4_i \lm_i^{-1}a_1(u_i,\f).
$$
Using (10) we obtain
$$
 |(u_i,\f)_{L_2(0,b)}|\le\e^2_i a_1^{1/2}(\f,\f).
$$
It follows immediately that $u_i\to 0$ in $L_2(0,b)$ weakly.

Throughout the proof, $\mathcal H_a$ denotes the space
$\{v\in H^2(a,0): v(a)=v'(a)=0\}$.
Passing to the limit as $\e_i\to 0$ in (5) we obtain
\begin{equation}
 a_0(v,\f)-\lm(v,\f)_{L_2(a,0)}=0,\quad \f\in H^2_0(a,b).
 \label{tag11}
\end{equation}
Since the restriction $\f|_{(0,a)}$ is an element of $\mathcal H_a$,
the identity (11) corresponds to the eigenvalue problem (9). This finishes the proof,
because $v$ is a nonzero function.
\hfill $\Box$

\vspace{5mm}

Consequently, the limit frequencies $\om$ are generated by means of stiffer part of a vibrating system.

\section{Asymptotic expansions of high frequency vibrations. The leading terms}
By $\om_\e$ we denote an eigenfrequency $\lm_\e^{1/4}$ of vibrating system
(1)--(4).  We search for $\omega_\e$ with the asymptotic expansion
\begin{equation}
 \om_\e \sim\sum\limits_{k=0}^{\infty}\e^k\om_k,
 \qquad \omega_0\not =0.
 \label{tag12}
\end{equation}
Let us consider a one-parameter family of the following vector-functions
$$
N(\kappa,x)=\left(\cos \kappa S(x),\, \sin \kappa S(x),\,  e^{-\kappa S(x)},
\, e^{\kappa (S(x)-S(b))}\right), \qquad
x\in (0,b),\quad \kappa\in \mathbb R,
$$
where
$$
S(x)=\omega_0 \int_0^xk_1^{-1/4}(t)\,dt.
$$
Let $<\cdot,\cdot>$ be an inner product in $\mathbb R^4$.
We postulate the expansions of $u_{\e}$ in the form
\begin{equation}
 u_\e(x) \sim
 \left\{ 
 \begin{array}{ll}
 \sum\limits_{k=0}^{\infty}\e^kv_k(x),&\quad x\in (a,0),\\
        \sum\limits_{k=0}^{\infty}
            \e^k<f_k(x), N\left(\e^{-1},x\right)>,
            &\quad x\in (0,b),
 \end{array}
 \right.
 \label{tag13}
\end{equation}
with $f_k: (0,b)\to\mathbb R^4$.

\begin{rem}
On the interval $(0,b)$ we have the equation
$$
  \e^4\frac{d^2}{dx^2}\left(k_1\frac{d^2u_\e}{dx^2}\right)-
   \om_\e^4u_\e=0
$$
with a small parameter at the highest order derivative. 
Thus the construction of expansions (13) 
on the interval $(0,b)$ was motivated by [13]
with
the method of WKB-approximations, which are also known as
short-wave approximations.
Hence, we can find a solution $u_\e$  in the form
$$
  e^{-\frac{S(x)}{\e}}\sum_{k=0}^{\infty}\e^ka_k(x),
$$
where  function $S(x)$ has to satisfy the eikonal equation
$
 k_1(x)S'(x)^4-\om_0^4=0.
$
Since there exist 4 different solutions of the eikonal equation with respect to
$S'$, we have introduced  the function $N(\kappa,\cdot)$
combining in it's components a fundamental set of solutions.
\end{rem}

Substituting (13) for $u_\e$ in conditions (4), we obtain
\begin{align}
 k_0(0)\, v''_0(0)+\dots &\sim \e^2\, k_1(0)\, S'(0)^2 <f_0(0),
   T^2N\left({\e}^{-1},0\right)>+\dots \,,
   \label{tag14}\\
 (k_0v''_0)'(0)+\dots &\sim \e\, k_1(0)\, S'(0)^3 <f_0(0),
   T^3N\left({\e}^{-1},0\right)>+\dots\,\,.
   \label{tag15}
\end{align}
Above we take into account the equality
$$
  \frac{d}{dx}N\left({\e}^{-1},x\right)=
     {\e}^{-1}{S'(x)}\, TN\left({\e}^{-1},x\right)
$$
with an orthonormal matrix
$$
  T=\left(\begin{array}{cc}
            {T_1}&{0}\\
            {0}&{T_2}
        \end{array}
  \right),\quad\text{where}\quad
  T_1=\left(\begin{array}{cc}
            {0}&{-1}\\
            {1}&{0}
        \end{array}
  \right),\quad
  T_2=\left(\begin{array}{cc}
            {-1}&{0}\\
            {0}&{1}
        \end{array}
  \right).
$$
By (14) and (15), we obtain  $v''_0(0)=0$ â  $v'''_0(0)=0$.
Substituting expansions (12), (13) for $\om_\e$ and $u_\e$ in equation (1)
and boundary value conditions (2), we have the eigenvalue problem
\begin{equation}
\begin{array}{c}
\displaystyle
 \frac{d^2}{dx^2}\left( k_0(x)\frac{d^2v_0}{dx^2}\right)-\om_0^4\,v_0(x)=0,
 \quad x\in (a,0),
 \\[2mm]
\displaystyle
 v_0(a)=v'_0(a)=0,\quad v''_0(0)=v'''(0)=0.
 \end{array}
\label{tag16}
\end{equation}
Let $\om_0$ and $v_0$ be an eigenfrequency and an eigenfunction of the problem
(compare with Th.1). Note that every eigenvalue $\om_0$ is simple.
Suppose that the function $v_0$ satisfies condition $\|v_0\|_{L_2(a,0)}=1$.

In order to find the leading term $f_0$ of expansion (13), we substitute
the second of series (13) in equation (1).  In particular, we obtain
\begin{eqnarray} 
&
(k_1{S'}^4-\om_0^4)f_0=0,
&
 \label{tag17}
 \\
 &
 4k_1{S'}^3Tf'_0+(2k_1'{S'}^3T+6k_1{S'}^2S''T
 -4\om^3_0\om_1 I)f_0=-(k_1{S'}^4-\om_0^4)f_1.
 &
 \label{tag18}
\end{eqnarray}

Since the function $S$ satisfies the eikonal equation $k_1{S'}^4-\om_0^4=0$
then (17) holds. Hence (18) is a homogeneous system of linear differential equations with respect to the function $f_0$.

Taking into account conditions (3) and boundary value conditions (2) at $x=b$,
we can write
\begin{eqnarray} 
&
f'_0=A(x)f_0,\quad x\in (0,b),
&
 \nonumber
 \\
 &
 \langle f_0(0), \,N(\e^{-1}, 0) \rangle = v_0(0),\quad
 \langle f_0(0), \,TN(\e^{-1}, 0) \rangle = 0,
 &
\label{tag19}
 \\
 &
 \langle f_0(b), \,N(\e^{-1}, b) \rangle = 0,\quad
 \langle f_0(b), \,TN(\e^{-1}, b) \rangle = 0.
 &
 \nonumber
\end{eqnarray}
with a matrix
$$
 A=\left(\begin{array}{cccc}
        {-\frac{k'_1}{8k_1}}&{\frac{\om_1}{\root{4}\of{k_1}}}&{0}&{0}\\
        {-\frac{\om_1}{\root{4}\of{k_1}}}&{-\frac{k'_1}{8k_1}}&{0}&{0}\\
        {0}&{0}&{-\frac{k'_1}{8k_1}-\frac{\om_1}{\root{4}\of{k_1}}}&{0}\\
        {0}&{0}&{0}&{-\frac{k'_1}{8k_1}+\frac{\om_1}{\root{4}\of{k_1}}}
        \end{array}
  \right).
$$
Note that the matrix $A$ depends on a parameter $\om_1$. We shall define $\om_1$ below.

On the one hand  problem (19) depends on $\e$ by means of
the boundary value conditions, on the other hand 
that one is an ill-posed problem.
To resolve both these problems, we consider a discrete set
of small parameter $\e$. We shall choose below a specific sequence $\e_p\to 0$
of a small parameter and asymptotic expansions will have a discrete character
with respect to $\e$.
That agrees with a discrete phenomenon of high frequency vibrations.

\begin{lemma} Let $w:(0, b)\to\mathbb R^4$ be a
smooth vector-function and $\sigma$ be a vector in $\mathbb R^4$.
There exists a infinitely small sequence $\{\e_p\}_{p=1}^\infty$
such that the problem
\begin{eqnarray} 
&
 y'(\e,x)=A(x)y(\e,x)+w(x),\quad x\in (0,b),
 &
 \label{tag20}
 \\
 &
 \langle y(\e,0),\,N(\e^{-1}, 0) \rangle  =\sigma_1,\quad
 \langle y(\e,0),\,TN(\e^{-1}, 0) \rangle =\sigma_2,
 &
 \label{tag21}
\\ 
&
\langle y(\e,b),\,N(\e^{-1}, b) \rangle =\sigma_3,\quad
 \langle y(\e,b),TN(\e^{-1}, b) \rangle =\sigma_4,
 &
 \label{tag22}
\end{eqnarray}
has a unique solution $y(\e_p,\cdot)$  for $p\in \mathbb N$.
The family of solutions $\{y(\e_p,\cdot)\}_{p=1}^\infty$ holds the inequality
$$
\| y(\e_p,\cdot)-y_*\|_{C^1}\le C e^{-\frac M{\e_p}}
$$
for a  smooth function $y_*:[0,b]\to \mathbb R^4$.
The constants $C$ and  $M$ are independent of 
parameter $\e$.
\end{lemma}

\proof
The fundamental matrix of (20) is
$$
  \Phi(x)=k_1^{-1/8}(x)\left(\begin{array}{cccc}
         {\cos\frac{\om_1}{\om_0}\displaystyle S(x)}&{\sin\frac{\om_1}{\om_0}\displaystyle S(x)}&{0}&{0}\\
         {-\sin\frac{\om_1}{\om_0}\displaystyle S(x)}&{\cos\frac{\om_1}{\om_0}\displaystyle S(x)}&{0}&{0}\\
         {0}&{0}&{e^{-\frac{\om_1}{\om_0}S(x)}}&{0}\\
         {0}&{0}&{0}&{e^{\frac{\om_1}{\om_0}(S(x)-S(b))}}
\end{array}\right).
$$
Therefore we have a representation of the general solution
$$
y(x)=\Phi(x)(\beta+h(x)),
$$
where $\beta$ is a constant vector
and  $h(x)=\int_0^x\Phi^{-1}(t)w(t)\,dt$. %
Suppose a vector-function $y(\e,x)=\Phi(x)(\beta_\e+h(x))$ is a solution of (20)--(22).
Seting $\gamma_\e=\e^{-1}+\om_1\om_0^{-1}$, it is easy to check that
\begin{equation}
 \Phi^t(x)N(\e^{-1}, x)= k_1^{-1/8}(x) N(\gamma_\e, x),
 \label{tag23}
\end{equation}
where $\Phi^t$ is a transposed matrix.
Then
$$
 <y(\e,x),N(\e^{-1},x)>=
 k_1^{-1/8}(x)<\beta_\e+h(x),N(\gamma_\e, x)>,
$$
and we can write (21), (22) in the form
\begin{equation}
\begin{array}{rl}
 \langle \beta_\e, \kern0.9em N(\gamma_\e,0) \rangle &=m_0\sigma_1,
 \\
 \langle \beta_\e,\,TN(\gamma_\e,0) \rangle &=m_0\sigma_2,
 \\
 \langle \beta_\e, \kern1em N(\gamma_\e,b) \rangle &=m_1\sigma_3-<h(b),\,N(\gamma_\e,b)) \rangle, 
 \\
 \langle \beta_\e,\, TN(\gamma_\e,b) \rangle &
 = m_1 \sigma_4 - \langle h(b),\,TN(\gamma_\e,b)) \rangle ,
\end{array}
\label{tag24}
\end{equation}
where $m_0=k_1^{1/8}(0)$ and $m_1=k_1^{1/8}(b)$.
Note that the matrix $\Phi^t$ commutes with $T$ and, moreover, $h(0)=0$.

Hence, the vector $\beta_\e$ is a solution of a linear algebraic
system with a matrix
$$
 G(\gamma_\e)=\left(\begin{array}{cccc}
           {1}&{0}&{1}&{e^{-\gamma_\e S(b)}}\\
               {0}&{1}&{-1}&{e^{-\gamma_\e S(b)}}\\
               {\cos\gamma_\e S(b)}&{\sin\gamma_\e S(b)}&{e^{-\gamma_\e S(b)}}&{1}\\
               {-\sin\gamma_\e S(b)}&{\cos\gamma_\e S(b)}&{-e^{-\gamma_\e S(b)}}&{1}
              \end{array}\right).
$$

Since $S(b)\not= 0$, the determinant
$$
\det G(\gamma_\e)=-2\cos\gamma_\e S(b)+2e^{-\gamma_\e S(b)}
\bigl(2-e^{-\gamma_\e S(b)}\cos\gamma_\e S(b)\bigr)
$$
does not vanish for all $\e>0$.

Let us fix $\delta\in [0, 2\pi)$ and choose a sequence $\e_p$ from
the set of conditions
$\gamma_{\e_p} S(b) = \delta+2\pi p$ for  $p=1,2,\dots \,\,$.
Namely, let
\begin{equation}
 \e_p=\frac{\om_0 S(b)} {\om_0(\delta+2\pi p)-\omega_1 S(b)}
 \label{tag25}
\end{equation}
for all $p\ge p_0$. Here $p_0$ is the smallest natural number
such that the denominator of (25) is positive. Without restriction of generality we set $p_0=1$.
We denote $\gamma_p=\gamma_{\e_p}$. Since $\gamma_p\to\infty$ as $\e_p\to0$,
the matrix
$G(\gamma_p)$ is an exponentially small perturbation of a matrix
$$
 G_0=\left(\begin{array}{cccc}
           {1}&{0}&{1}&{0}\\
               {0}&{1}&{-1}&{0}\\
               {\cos\delta}&{\sin\delta}&{0}&{1}\\
               {-\sin\delta}&{\cos\delta}&{0}&{1}
              \end{array}\right).
$$

Since
$N(\gamma_p,b)=\left(\cos\delta, \sin\delta, e^{-\gamma_p S(b)}, 1 \right)$
then
the right-hand side of (24) is an exponentially small perturbation of vector
$$
g = (m_0\sigma_1,\, m_0\sigma_2,\, m_1\sigma_3-
\langle h(b),\,N_{\delta}\rangle ,
\, m_1\sigma_4-\langle h(b),\, TN_{\delta}\rangle ),
$$
where  $N_{\delta}$ differs from vector $N(\gamma_p,b)$
by the third component only.

Let us suppose that $\delta$ is not equal to ${\pi}/{2}$ and $3\pi/2$. Then the matrix $G_0$ is non-degenerate.
From the theory of finite-dimensional perturbations we obtain
$$
 \|\beta_{\e_p}-\beta_*\|_{\mathbb R^4}\le C e^{-\gamma_p S(b)},
$$
where  $\beta_*$  is a solution of $G_0\beta=g$.
Let $y_*(x)=\Phi(x)(\beta_*+h(x))$, then
\begin{equation}
\| y(\e_p,\cdot)-y_*\|_{C^1}\le \|\Phi\|_{C^1}
\|\beta_{\e_p}-\beta_*\|_{\mathbb R^4},
\label{tag26}
\end{equation}
where $\|\Phi\|_{C^1}=\max\limits_{x\in(0,b)}(\|\Phi(x)\|+\|\Phi'(x)\|)$.
Note that $\gamma_\e\ge c_o\e^{-1}$ with a positive constant
$c_0$, which proves the Lemma.
\hfill $\Box$

\begin{rem} 
From now on, we say that  $y_*$
is a solution of (20)--(22) with neglect of exponentially small terms.
However the choice of a sequence $\e_p$ is non-unique and depends on
$\delta$ and  $\omega_1$. We shall define $\omega_1$ at the next step
but, on the other hand, we shall keep the  dependence on $\delta$,
which will be a deformation parameter.
The approximation of the limit function $v_0$ by eigenfunctions
$u^\e_{n(\e)}$ is ambiguously determined. Hence we can not define
$\delta$ uniquely. Consequently, for each  $\delta$  we  construct
the expansions of $u^\e$, though they are asymptotically equivalent
as $\e\to 0$.
\end{rem}

Let us return to the study of problem (19), that is a sub-case of (20)--(22)
with the right-hand side $w=0$ and $\sigma=(v_0(0),0,0,0)$.
Since (19) is a homogeneous system, we obtain $f_0=\Phi\beta_0$, where
$
 \beta_0=\frac{1}{2}\displaystyle k_1^{1/8}(0)v_0(0)
 \left(1-\tg\delta,1+\tg\delta,1+\tg\delta,-\frac{1}{\cos\delta}\right)
$
is a solution of the corresponding system with matrix $G_0$.

In order to calculate the first-order correction $\om_1$,
we consider the problem for $v_1$
\begin{equation}
\begin{array}{rcl} 
&\displaystyle
\frac{d^2}{dx^2}\left( k_0(x)\frac{d^2v_1}{dx^2}\right)-\om_0^4v_1=
   4\om_1\om_0^3v_0,\quad x\in (a,0),
   &
   \\
   &\displaystyle
 v_1(a)=v'_1(a)=0,\quad v''_1(0)=0,\quad
   v'''_1(0)=-\frac{k_1^{1/4}(0)}{k_0(0)}\om_0^3v_0(0)(1+\tg\delta).
   &
\end{array}
\label{tag27}
\end{equation}
Since $\om_0^4$ is a simple eigenvalue of (16), there exists a solution
of (27) under the condition
$
 \om_1=-\frac{1}{4}k_1^{1/4}(0)v_0^2(0)(1+\tg\delta).
$
Let us choose a solution $v_1$ such that $(v_1,v_0)_{L_2(a,0)}=0$.

Now by (25), it follows that
\begin{equation}
 \e_p(\delta)=\frac{4\om_0 S(b)}
  {4\om_0(\delta+2\pi p)+k_1^{1/4}(0)v_0^2(0)S(b)(1+\tg\delta)},
  \label{tag28}
\end{equation}
where  $\delta\in [0,2\pi)\setminus \{\pi/2,3\pi/2\}$.
Let us denote by $\mathcal E_{\delta}$
 the sequence (28) for fixed  $\delta$. 
 We shall construct the asymptotic expansions only for
$\e\in \mathcal E_{\delta}$.
Hence, we have just found  $\om_0$, $v_0$, $f_0$, tgcorrections
$\om_1$, $v_1$ and the sequence $\mathcal E_{\delta}$.
Recall that $f_0$, $\om_1$ and $v_1$ depend on $\delta$.

\section{Complete asymptotics of high frequency vibrations}
Let us find the general terms  $f_k$, $\om_{k+1}$ and $v_{k+1}$
of expansions (12), (13) for $k\geq 1$.
In the same way as for the vector $f_0$, we obtain the boundary value
problem for $f_k$. With neglect of the
exponentially small terms in the boundary conditions, we can write
the  problem  in the form
\begin{equation}
\begin{array}{rl}
 f'_k &=A(x)f_k+w_k,\quad x\in (0,b),
 \\
 \langle f_k(0),\,N(\e^{-1},0) \rangle  &=v_k(0),\qquad 
 \langle f_k(b),\,N(\e^{-1},b) \rangle =0,
\label{tag29}
 \\
 \langle f_k(0),\,TN(\e^{-1},0) \rangle  &=\om_0^{-1}(k_1^{1/4}(0)v'_{k-1}(0)-
 k_1^{1/8}(0)\langle \Phi^{-1}(0)f'_{k-1}(0),N_0 \rangle ),
 \\
 \langle f_k(b),\,TN(\e^{-1},b) \rangle  &=
 -\om_0^{-1}k_1^{1/8}(b)\langle \Phi^{-1}(b)f'_{k-1}(b),N_{\delta} \rangle,
\end{array}
\end{equation}
where  $N_0=(1,0,1,0)$
and the $N_{\delta}$ is defined in the proof of Lemma 3.

For each $k$ the right-hand side $w_k$ of (29) is a smooth function, namely,
\begin{equation*}
\gathered
 w_k=\frac{1}{4k_1{S'}^3}\Bigl(
 \sum_{m=0}^{k-1}\lm_{k-m}T^3f_m-
 S'T^2P_k-\frac{d}{dx}(k'_1(S')^2Tf_{k-1}+
 \\
 3k_1S'S''Tf_{k-1}
 +3k_1(S')^2Tf'_{k-1}+T^3P_{k-1})\Bigr),
\endgathered
\end{equation*}
where
$
 P_l=k_1S'T^3f''_{l-1}+\frac{d}{dx}(2k_1S'T^3f'_{l-1}+k_1S''T^3f_{l-1}+
 k_1f''_{l-2}).
$
The proof is by induction on $k$.
According to Lemma 3  there exists a solution of (29) for
$\e\in \mathcal E_{\delta}$.

Now we can find $\om_{k+1}$ and $v_{k+1}$:
\begin{equation}
\begin{array}{c}
\displaystyle
\frac{d^2}{dx^2}\left( k_0(x)\frac{d^2v_{k+1}}{dx^2}\right)-\om_0^4v_{k+1}=
   \sum_{m=1}^{k+1} \lm_{m} v_{k-m+1},\quad x\in (a,0),
   \\
   v_{k+1}(a)=0,\qquad v'_{k+1}(a)=0,
   \label{tag30}
   \\
   (k_0v''_{k+1})(0)=
   \langle \Phi^{-1}(\om_0^2k_1^{3/8}T^2f_{k-1}+
   k_1^{-1/8}Q_{k-2}),N_0 \rangle  \Big|_{x=0},
   \\
   (k_0v''_{k+1})'(0)= \langle \Phi^{-1}(\om_0^3k_1^{1/8}Tf_k+
   k_1^{-1/8}R_{k-1}),N_0\rangle \Big|_{x=0},
\end{array}
\end{equation}
where
\begin{align*}
 Q_k &=k_1\frac{d}{dx}(S'T^3f_k+f'_{k-1})+S'T^3f'_k,
 \\
 R_k &=S'T^3Q_{k-1}+\frac{d}{dx}(\om_0^2k_1^{1/2}T^2f_k+Q_{k-1})
\end{align*}
and
$
 \lm_{m}=\sum\om_i\om_j\om_l\om_s
$
with $i+j+l+s=m$.
Since $\om_0^4$ is an eigenvalue of (16), boundary value problem
(30) has no solution for an arbitrary right-hand side. We can write
the existence condition in the form
\begin{align*} 
\om_{k+1}=\frac{1}{4\om_0^3}<v_0\Phi^{-1}
( &
 \om_0^2 k_1^{3/8}T^2f_{k-1}+
   k_1^{-1/8}Q_{k-2})-
   \\
 &  
   v'_0\Phi^{-1}(\om_0^3k_1^{1/8}Tf_k+
   k_1^{-1/8}R_{k-1}),N_0> \Big|_{x=0}.
\end{align*}
Let us choose a solution  $v_{k+1}$ of (30) such that $(v_{k+1},v_0)_{L_2(a,0)}=0$.

Hence,  the algorithm scheme of asymptotics (12),~(13) is
$$
 \om_0\to v_0\to
 f_0\to\om_1\to\mathcal E_{\delta}\to v_1
 \to\dots\to
 f_k\to\om_{k+1}\to v_{k+1}\to\dots \,\, .
$$
Note that all terms, except for $\om_0$ and $v_0$,
depend on parameter $\delta$.

\section{Justification of asymptotics}
The nonstandard object of investigation, namely, the sequences of
eigenfunctions $u^\e_{i(\e)}$, changes the classical scheme of justification.
Note that in Sections 4 and 5 we have constructed 
series (12), (13), however we have not
defined the object that is approximated by ones yet.

Only by using the formal series (12) we shall define a sequence $u^\e_{k(\e)}$
that supports the high frequency vibrations with the limit frequency $\om_0$
and the limit form $v_0$.

Let us fix  $n\in\mathbb N$ and introduce a sequence of real numbers  $\{\lm_\e^{(n)}\}_{\e\in\mathcal E_{\delta}}$
and sequence of functions  $\{u_\e^{(n)}\}_{\e\in\mathcal E_{\delta}}$ in
$H^2(a,b)$. Namely, for each $\e\in \mathcal E_{\delta}$ we set
\begin{align}
 \lm^{(n)}_\e & = \lm_0 + \e \lm_1 + \dots \e^n \lm_n ,\qquad
 \lm_s=\kern-4pt\sum_{i+j+k+l=s}\kern-6pt\om_i\om_j\om_k\om_l,
 \label{tag31}\\
 u^{(n)}_\e(x)
 &=\left\{\begin{array}{ll}
           v_0(x)+\e v_1(x)+ \dots + \e^n v_n(x), &\quad x\in (a,0),\\
           \sum\limits_{i=0}^n\e^i<f_i(x),
               N(\e^{-1},x)>,& \quad x\in (0,b),
           \end{array}
           \right.
 \label{tag32}
\end{align}
where the numbers $\om_k$, the functions $v_k$, the vectors $f_k$ and the set
$\mathcal E_{\delta}$ are defined in Sections 4 and 5.

\begin{lemma}
There exists a sequence $\{\lm^\e\}_{\e\in\mathcal E_{\delta}}$ of eigenvalues
for (1)--(4) such that
\begin{equation}
 |\lm^{(n)}_\e-\lm^{\e}|\leq C_n\e^{n+1},\quad \e\in\mathcal E_{\delta}.
 \label{tag33}
\end{equation}
\end{lemma}

\proof
Denoting by $\mathcal H_\e$ the space $H^2_0(a,b)$ with the
inner product $a_\e(\cdot,\cdot)$,
let us introduce an operator $A_\e$ such that
$$
 a_\e (A_\e u,v)=(u,v)_{L_2(a,b)},\qquad v\in\mathcal H_\e.
$$
Note that $A_\e$ is self-adjoint and compact for all $\e>0$.
Then we can write problem (1)--(4) in the form
$$
  A_\e u_\e-(\lm^\e)^{-1} u_\e=0.
$$
Substituting sequences (31), (32) in problem (1)--(4), we obtain
\begin{equation}
\begin{array}{c}
\displaystyle
\frac{d^2}{dx^2}\left(k_\e(x)\frac{d^2}{dx^2}u^{(n)}_\e\right)-
 \lm^{(n)}_\e u^{(n)}_\e
 =F_n(\e,x),
 \quad x\in (a,b),
 \\
\displaystyle
 u^{(n)}_\e(a)=0,
 \quad 
 \frac{d}{dx}u^{(n)}_\e(a)=0,
 \quad
 u^{(n)}_\e(b)=g^{(1)}_n(\e),
 \quad
 \frac{d}{dx}u^{(n)}_\e(b)
 =g^{(2)}_n(\e),
 \\[1mm]
\displaystyle
 \left[ u^{(n)}_\e \right]_0 = g^{(3)}_n(\e),
 \quad
 \left[ \frac{d}{dx}u^{(n)}_\e \right]_0 = h_n(\e),
 \\[1mm]
\displaystyle
 \left[k_\e\frac{d^2}{dx^2}u^{(n)}_\e\right]_0
 =z_n^{(1)}(\e),
 \quad
 \left[\frac{d}{dx}(k_\e\frac{d^2}{dx^2}u^{(n)}_\e)\right]_0
 =z_n^{(2)}(\e),
\end{array}
 \label{tag34}
 \end{equation}
where  $[f]_0$ is a jump of a function $f$ at $x=0$.
The right-hand sides of problem (34) satisfy the inequalities
\begin{equation}
\begin{array}{c}
\displaystyle
\|F_n(\e,x)\|_{C^0(a,b)}\leq C_n\e^{n+1},
\quad
|g^{(i)}_n(\e)|\leq C_n e^{-\frac{M}\e},
\quad i=1,2,3,
\\[2mm]
\displaystyle
|h_n(\e)|\leq C_n\e^n,   
\quad
|z_n^{(i)}(\e)|\leq C_n\e^{n+1},
\quad i=1,2.
\end{array}
\label{tag35}
\end{equation}
The function $u^{(n)}_\e$ does not belong to the space $\mathcal H_\e$, because it has a point of discontinuity at $x=0$.  
Taking into account (35) we can choose
a function $\f^{(n)}_\e$ 
such that $u^{(n)}_\e+\f^{(n)}_\e \in\mathcal H_\e$ and
\begin{equation}
  \max_{x\in (a,b)}\left(\left|\f^{(n)}_\e\right|+
    \left|\frac{d}{dx}\f^{(n)}_\e\right|+
    \left|k_\e^{1/2}\frac{d^2}{dx^2}\f_\e^{(n)}\right|\right)\leq C_n\e^n.
    \label{tag36}
\end{equation}
Let $V^{(n)}_\e=\kappa_\e(u^{(n)}_\e+\f^{(n)}_\e)$, where
$\kappa_\e$ is a normalizing constant such that
 $\|V^{(n)}_\e\|_{\mathcal H_\e}=1$. It is easy to check that
$\kappa_\e\ge\kappa_0>0$.
It follows from (35), (36) that
\begin{equation}
 \left\|\left(A_{\e}-(\lm_\e^{(n)})^{-1}I\right)
  V_\e^{(n)}\right\|_{\e}
 \leq K_n\e^n, \quad \e\in\mathcal E_{\delta},
 \label{tag37}
\end{equation}
where $K_n$ is a constant independent of $\e$.
Hence, according to the Vishik-Lusternik lemma [14] there exists the eigenvalue
$(\lm^\e)^{-1}$ of the operator  $A_{\e}$ such that
$$
 \left|\frac{1}{\lm^\e}-\frac{1}{\lm_\e^{(n)}}\right|\leq K_n\e^n.
$$
Applying this inequality for the value $n+1$ instead of $n$, we obtain
$$
 |\lm^{\e}-\lm_\e^{(n)}|\leq C_n\e^{n+1},\qquad
 \e\in\mathcal E_{\delta}.
$$
\hfill $\Box$

\begin{lemma}
There exists a positive number $d$ such that for each  $\e\in\mathcal E_{\delta}$ the interval
$I_\e=(\lambda_\e^{(n)}-d\e,\lambda_\e^{(n)}+d\e)$
contains exactly one eigenvalue $\lm_{k(\e)}^\e$ of problem (1)--(4).
Moreover, the eigenvalue number $k(\e)$ satisfies the inequality
\begin{equation}
 a_0\e^{-1}\le k(\e)\le a_1\e^{-1},
 \label{tag38}
\end{equation}
where $a_0$, $a_1$  are constants independent of $\e$.
\end{lemma}

\proof
According to Section 2, we have
$$
 \lm_m^\e=\e^4\mu_m(1+\alpha_m(\e)),\quad\e\to 0,
$$
where $\mu_m$ is an eigenvalue of (8) and
$\alpha_m(\e)=o(1)$ as  $\e\to 0$.
Moreover, the following asymptotics  hold
$$
 \mu_m=m^4(c_0+\tau(m))
 \quad
 \text{with}
 \quad \tau(m) = o(1)
 \quad
 \text{as}
 \quad
  m\to\infty.
$$
Then
\begin{equation}
 \lm_m^\e=(\e m)^4(c_0+\tau(m))(1+\alpha_m(\e)).
 \label{tag39}
\end{equation}
By Lemma 4, for each $\e\in\mathcal E_\delta$ there exists 
a number $k=k(\e)$
such that
$$
 \lm_0-b_0\e\le \lm_{k(\e)}^\e\le
 \lm_0+b_0\e
$$
with $b_0>0$. We conclude from (39)  that
$$
 \lm_0-b_0\e\le (\e k)^4(c_0+\tau(k))(1+\alpha_k(\e))\le
 \lm_0+b_0\e,
$$
and finally (38) holds.
According to (39) we have
$$
 \lm_{k(\e)}^\e=c_0\e^4k(\e)^4+o(\e^4),\qquad \mathcal E_\delta\ni\e\to 0.
$$
Then the distance between the two neighboring eigenvalues
that are the
closest to the point $\lm_0$ admits the estimate
\begin{equation}
 |\lm_{k(\e)+1}^\e-\lm_{k(\e)}^\e|=4c_0\e^4k(\e)^3+o(\e^4),
 \qquad \mathcal E_\delta\ni\e\to0.
 \label{tag40}
\end{equation}
Taking into account the asymptotic behaviour of $k(\e)$, 
we obtain the existence of the interval $I_\e$ of the length $2d\e$ such that exactly one eigenvalue $\lm^\e$ of problem (1)--(4) belongs to it. 
The case of the point $\lm_\e^{(n)}$ being
a midpoint of $I_\e$ is impossible since this contradicts (33) and (40).
\hfill $\Box$

\vspace{5mm}

Now we can introduce the object that is approximated by formal series (13)
constructed in Sections 4 and 5.
Let $\{\lm^\e_{k(\e)}\}_{\e\in\mathcal E_\delta}$
be the sequence which is defined in Lemma 5, where $\lm^\e_{k(\e)}$ is the nearest eigenvalue to the value $\lm_\e^{(n)}$. Let $\{u^\e_{k(\e)}\}_{\e\in\mathcal E_\delta}$ be
the sequence of the corresponding eigenfunctions,
$\|u^\e_{k(\e)}\|_{\mathcal H_\e}=1$.

\begin{theorem}
As $\mathcal E_\delta\ni\e\to0$
the sequence $\{ u^\e_{k(\e)}\}_{\e\in\mathcal E_\delta}$ supports
 the high frequency vibrations with the limit
frequency
$\om_0$  and the limit form $v_0$.  Moreover,
\begin{equation}  
\|u^\e_{k(\e)}-V^{(n)}_\e\|_\e\leq C_n\e^{n+1}, 
  \qquad\e\in\mathcal E_\delta,
 \label{tag41}
 \end{equation}
for $n=0,1,\dots\,\,$. In particular, the following inequalities hold
\begin{equation}
\begin{array}{c}
\displaystyle
  \left\|u^\e_{k(\e)}(x)-\kappa_\e{\sum}_{k=0}^n v_k(x)\e^k\right\|_{C^1(a,0)}
  \leq C_n\e^{n+1},
\\
\displaystyle
  \left\|u^\e_{k(\e)}(x)-\kappa_\e
  {\sum}_{k=0}^n\langle f_k(x),N(\e^{-1},x)\rangle\e^k\right\|_{C^1(0,b)}
  \leq C_n\e^{n-1},
\end{array} 
\label{tag42}
\end{equation}
where $\kappa_\e=\|V^{(n)}_\e\|_\e^{-1}$ is a normalizing constant and
$\kappa_\e\ge\kappa_0>0$.
\end{theorem}

\proof
By Lemma 7, we choose $d_0>0$ such that  $d_0\e$-neighborhood of the point
$(\lm_\e^{(n)})^{-1}$ contains exactly one eigenvalue
$(\lm_{k(\e)}^\e)^{-1}$ of the operator $A_\e$.
From inequality (37), we obtain (see [14])
\begin{equation}
 \|u_\e-V_\e^{(n)}\|_\e\leq K_n d_0^{-1}\e^{n-1},
 \label{tag43}
\end{equation}
where $u_\e$ is a normalized eigenfunction of $A_\e$ associated with
the eigenvalue $(\lm_{k(\e)}^\e)^{-1}$.
Hence, $u_\e=\pm u^\e_{k(\e)}$.
There is no loss of generality in assumption that (43) holds for
the function $u^\e_{k(\e)}$.
Taking into account that the terms of series (13) are bounded in
$\mathcal H_\e$-norm, from (43) for the number $n+2$, we obtain (41).

In addition, by (36) and the definition of norm, we obtain
$$
 \gathered
  \|u^\e_{k(\e)}-\kappa_\e u_\e^{(n)}\|_{H^2(a,0)}\leq C_n\e^{n+1}, \qquad
  \|u^\e_{k(\e)}-\kappa_\e u_\e^{(n)}\|_{H^2(0,b)}\leq C_n\e^{n-1}.
 \endgathered
$$
Estimate (42) follows from the Sobolev embedding theorem.

It is easy to check that the constant $\kappa_\e$ has a nonzero limit as
$\e\to0$. Then sequence $\{ u^\e_{k(\e)}\}_{\e\in\mathcal E_\delta}$ converges
to the eigenfunction $\kappa v_0$ of problem (8) in $(a,0)$, namely,
this sequence supports the high frequency vibrations.
\hfill $\Box$

\end{document}